\documentclass{amsart}
\usepackage{bm}

\newtheorem{theorem}{Theorem}[section]
\newtheorem{lemma}[theorem]{Lemma}
\newtheorem{proposition}[theorem]{Proposition}
\newtheorem{corollary}[theorem]{Corollary}
\theoremstyle{definition}
\newtheorem{definition}[theorem]{Definition}
\newtheorem{example}[theorem]{Example}

\theoremstyle{remark}
\newtheorem{remark}[theorem]{Remark}

\numberwithin{equation}{section}

\theoremstyle{plain}

\def\int{\mathop{\roman{int}}}

\def\1{^{-1}}

\def\Z{{\mathbf Z}}

\def\proof{{\bf Proof. }}
\def\endproof{\hfill \qed}

\numberwithin{equation}{section}
\begin{document}

\title[Group Actions and Covering Maps in the Uniform Category]
{Group Actions and Covering Maps in the Uniform Category}

\author{N.~Brodskiy}
\address{University of Tennessee, Knoxville, TN 37996}
%    Current address (if needed):
%\curraddr{}
\email{brodskiy@@math.utk.edu}
%\thanks{The first author was supported in part by NSF Grant \#000000.}

%    Information for second author (if needed):
\author{J.~Dydak}
\address{University of Tennessee, Knoxville, TN 37996}
\email{dydak@@math.utk.edu}
\thanks{The second-named author was partially supported
by the Center for Advanced Studies in Mathematics
at Ben Gurion University of the Negev (Beer-Sheva, Israel).}

\author{B.~LaBuz}
\address{University of Tennessee, Knoxville, TN 37996}
\email{labuz@@math.utk.edu}

%    Information for third author (if needed):
\author{A.~Mitra}
\address{University of South Florida, St. Petersburg}
\email{atish@stpt.usf.edu}


\date{February 26, 2008.}

\begin{abstract}
In \textit{Rips Complexes and Covers in the Uniform Category} \cite{BDLM} we define, following James \cite{Jam}, covering maps of uniform spaces and introduce the concept of generalized uniform covering maps. In this paper we investigate when these covering maps are induced by group actions.  Also, as an application of our results we present an exposition of Prajs' \cite{Pra} homogeneous curve that is path-connected but not locally connected.
\end{abstract}
\maketitle

\medskip
%Printed on \today.
\medskip
\tableofcontents

\section{Introduction}

\textit{Rips Complexes and Covers in the Uniform Category} \cite{BDLM} develops a theory of covering maps for uniform spaces. This theory includes two classes of maps.  Uniform covering maps, originally defined by James \cite{Jam}, are analogous to classical covering maps. In particular a uniform space has a universal uniform covering space if and only if it is uniform Poincare, i.e., is path connected, uniformly locally path connected, and uniformly semilocally simply connected. Generalized uniform covering maps were invented to deal with spaces that are not as nice. The concept of generalized paths is introduced as well as the related idea of a space being uniformly joinable.  It is then shown that a uniform space has a universal generalized uniform covering space if and only if it is uniformly joinable and chain connected.

This paper continues the development of the theory by investigating when uniform covering maps and generalized uniform covering maps arise from a group acting on a uniform space. Then, as an application, we provide an exposition of J.Prajs' \cite{Pra} example of a homogeneous curve that is path-connected but not locally connected. We thank Misha Levin for suggesting that we include such an exposition.

\section{Definitions and Results from ``Rips Complexes and Covers in the Uniform Category''}

We will discuss exclusively symmetric subsets $E$ of $X\times X$
(that means $(x,y)\in E$ implies $(y,x)\in E$) and the natural notation here (see \cite{Pla})
is to use $f(E)$ for the set of pairs $(f(x),f(y))$, where $f\colon X\to Y$ is a function.
Similarly, $f^{-1}(E)$ is the set of pairs $(x,y)$ so that $(f(x),f(y))\in E$ if $f\colon X\to Y$
and $E\subset Y\times Y$.
\par The {\bf ball $B(x,E)$ at $x$ of radius $E$} is the set of all $y\in X$
satisfying $(x,y)\in E$.
\par
A {\bf uniform structure} on $X$ is a family $\mathcal{E}$ of symmetric subsets $E$
of $X\times X$ (called {\bf entourages}) that contain the diagonal
of $X\times X$, form a filter (that means $E_1\cap E_2\in \mathcal{E}$
if $E_1,E_2\in \mathcal{E}$ and $F_1\in \mathcal{E}$ if $F_2\in \mathcal{E}$
and $F_2\subset F_1$),
and every $G_1\in \mathcal{E}$ admits $G\in \mathcal{E}$ so that
$G^2\subset G_1$ ($G^2$ consists of pairs $(x,z)\in X\times X$
so that there is $y\in X$ satisfying $(x,y)\in G$ and $(y,z)\in G$).
A {\bf base} $\mathcal{F}$ of a uniform structure $\mathcal{E}$
is a subfamily $\mathcal{F}$ of $\mathcal{E}$ so that for every entourage $E$ there is a subset
$F\in \mathcal{F}$ of $E$.

There are various equivalent definitions of a uniform covering map. In order to treat these, we need some preliminary definitions.

\begin{definition}
A map $f:X\to Y$ \textbf{generates the uniform structure on $Y$} if for each entourage $E$ of $X$, $f(E)$ is an entourage of $Y$
\end{definition}

Notice that if $f$ satisfies the above condition then it indeed generates the uniform structure on $Y$ in the sense that $\{f(E):E \textrm{ is an entourage of } Y\}$ is a basis for the uniform structure on $Y$. Also if $f$ generates the uniform structure on $Y$ then it is necessarily uniformly continuous and surjective.

\begin{definition}
Given an entourage $E$ of a uniform space $X$ and a subset $A\subset X$, $A$ is $\bm{E}$\textbf{-bounded} if for each $x,y\in A$, $(x,y)\in E$. Equivalently, $A\times A\subset E$.
\end{definition}

\begin{definition}
Given an entourage $E$ of $X$ the \textbf{Rips complex} $\bm{R(X,E)}$ is the subcomplex of the full complex over $X$ whose simplices are finite $E$-bounded subsets of $X$.
\end{definition}

We consider homotopy classes of paths in $R(X,E)$ joining two of
its vertices. Since the identity function $K_w\to K_m$, $K$ a
simplicial complex, from $K$ equipped with the CW (weak) topology
to $K$ equipped with the metric topology is a homotopy equivalence
(see \cite[page 302]{MarSeg}), it does not really matter which
topology we choose for $R(X,E)$.

The simplest path in $R(X,E)$ is the {\bf edge-path} $e(x,y)$ starting from
$x$ and ending at $y$ so that $(x,y)\in E$. Any path in $R(X,E)$ joining two vertices $x$ and $y$ can be realized, up to homotopy (see \cite[Section 3.4]{Spa}), as a
concatenation of edge-paths. Thus, each path in $R(X,E)$ can be realized by an $\bm{E}$\textbf{-chain} $x=x_1$, \ldots, $x_n=y$ such that $(x_i,x_{i+1})\in E$ for all $i<n$. Two paths in $R(X,E)$ represented by different $E$-chains with the same end-points are
homotopic rel. end-points if and only if one can move from one
chain to the other by simplicial homotopies: a new vertex $v$ can
be added or removed from a chain if and only if $v$ forms a
simplex in $R(X,E)$ with adjacent links of a chain (see
\cite[Section 3.6]{Spa}).

Given $f\colon X\to Y$ and an entourage $E$ of $X$ notice it induces a natural simplicial map $f_E\colon R(X,E)\to R(Y,f(E))$
by the formula $f_E(\sum\limits_{i=1}^n t_i\cdot x_i)=\sum\limits_{i=1}^n t_i\cdot f(x_i)$.

\begin{definition}\label{ChainLiftingProperty}
A surjective function $f\colon X\to Y$ from a uniform space $X$
has the \textbf{chain lifting property} if for any entourage $E$ of $X$ there is an entourage $F$ of $X$ such that any $f(F)$-chain in $Y$ starting from $f(x_0)$ can be lifted to an $E$-chain starting from $x_0$. 

A function $f\colon X\to Y$ from a uniform space $X$ has the \textbf{uniqueness of chain lifts property} if for every entourage $E$ of $X$ there is an entourage $F$ of $X$ such that any two $F$-chains $\alpha$ and $\beta$ satisfying $f(\alpha)=f(\beta)$ are equal if they originate from the same point.

We say that a function has the \textbf{unique chain lifting property} if it has the  chain lifting property and the uniqueness of chain lifts property.
\end{definition}

If a function $f:X\to Y$ has chain lifting then it generates a uniform structure on $Y$.

\begin{theorem} \label{EquivalentVersionsOfCovers}
Given a map $f:X\to Y$ that generates the uniform structure on $Y$, the following are equivalent.
\begin{itemize}

\item[1.] There is a basis $\mathcal{E}$ for the uniform structure on $X$ such that for each $E\in \mathcal{E}$, $f$ maps $B(x,E)$ bijectively onto $B(f(x),f(E))$.

\item[2.] There is a basis $\mathcal{E}$ for the uniform structure on $X$ such that for each $E\in \mathcal{E}$, the induced map $f_E\colon R(X,E)\to R(Y,f(E))$
is a simplicial covering map.

\item[3.] 
\begin{itemize}
\item[a.] For each entourage $E$ of $X$ there is an entourage $F$ of $Y$ such that if $(x,y)\in F$ and $f(x')=x$ there is a $y'\in X$ with $(x',y')\in E$
\newline and
\item[b.] There is an entourage $E_0$ of $X$ such that if $(x,y)\in E_0$ and $f(x)=f(y)$ then $x=y$.
\end{itemize}

\item[4.] $f$ has the unique chain lifting property.
\end{itemize}
\end{theorem} 

We call a simplicial map a \textbf{simplicial covering map} if it is a topological cover.

For the equivalence of (3) and (4), we have in particular that (3a) is equivalent to chain lifting and (3b) is equivalent to uniqueness of chain lifts. The entourage in (3b) is said to be \textbf{transverse} \cite{Jam} to $f$. 

\begin{definition}
A \textbf{uniform covering map} is a map $f:X\to Y$ between uniform spaces that generates the uniform structure on $Y$ and satisfies one of the equivalent conditions in \ref{EquivalentVersionsOfCovers}.
\end{definition}

A uniform covering map is a topological covering map and a uniform space has a universal uniform covering space if and only if it is a uniform Poincare space, i.e., is path connected, uniformly locally path connected, and uniformly semilocally simply connected. See \cite{BDLM} for the appropriate definitions.

Referring back to \ref{EquivalentVersionsOfCovers}, (1) helps to indicate that it is related to classical covering maps but is not very useful in practice. Statement (2) also gives a connection to classical covering maps and is interesting in its own right as Rips complexes are also vital for generalized uniform covering maps.  The characterization that is most useful in practice is (3).   Version (4) concisely states the essence of uniform covering maps and also relates uniform covering maps to generalized uniform covering maps.

The idea behind generalized uniform covers is to find a covering space theory that is useful for spaces that may not be path connected or locally path connected.  We shift the focus from paths to generalized paths. A generalized path is an element of the inverse limit of the Rips complexes $R(X,E)$ but can also be described in a nice way. A \textbf{generalized path} is a collection $\{[c_E]\}_E$ of homotopy classes of paths $[c_E]$ in $R(X,E)$ joining fixed $x\in X$ to $y\in X$ such that for all entourages $F\subset E$, $c_F$ is homotopic to $c_E$ in $R(X,E)$ rel. end-points. The space of generalized paths is denoted $GP(X)$. We normally restrict our attention to the pointed space $GP(X,x_0)$ consisting of generalized paths originating at a point $x_0\in X$.

To describe the uniform structure that we give $GP(X,x_0)$ let us first introduce some terminology. Let us say that two $E$-chains $c_E$ and $d_E$ starting at the same point are $\bm{E}$\textbf{-homotopic} if the concatenation $c_E^{-1}*d_E$ is $E$-homotopic to the edge path between their endpoints. Note that we necessarily have that their endpoints are $E$-close. Similarly let us say that two generalized paths $c=\{[c_E]\}_E$ and 
$d=\{[d_E]\}_E$ are $\bm{E}$\textbf{-homotopic} if $c_E$ and $d_E$ are $E$-homotopic. A generalized path is called $\bm{E}$\textbf{-short} if its $E$ term is $E$-homotopic rel. endpoints to the edge path in $R(X,E)$. Now, given an entourage $E$ of $X$ define a basic entourage $E^*$ of $GP(X,x_0)$ to be all pairs of generalized paths $(\alpha,\beta)$ that are $E$-homotopic, i.e., that have $\alpha^{-1}*\beta$ $E$-short.

We are now ready to define a generalized uniform covering map.

\begin{definition}
A map $f:X\to Y$ between uniform spaces that generates the uniform structure on $Y$ is a \textbf{generalized uniform covering map} if 
\begin{itemize}
\item[1.] $f$ has chain lifting
\item[2.] $f$ has approximate uniqueness of chain lifts
\item[3.] $f$ has generalized path lifting
\end{itemize}

A map $f:X\to Y$ has \textbf{approximate uniqueness of chain lifts} if for each entourage $E$ of $X$ there is an entourage $F$ of $X$ such that any two $F$-chains in $X$ originating at the same point are $E$-close if their images under $f$ are identical.  Two chains $c_1,\ldots, c_n$ and $d_1,\ldots, d_n$ are $\bm{E}$\textbf{-close} if $(c_i,d_i)\in E$ for each $i\leq n$.
\end{definition}

In this paper we do not deal with generalized path lifting.  Instead we use the following result.

\begin{proposition} \label{GenUnifCovForCompleteFibers}
If a map $f:X\to Y$ that generates the uniform structure on $Y$ has chain lifting, approximate uniqueness of chain lifts, and complete fibers then it is a generalized uniform covering map.
\end{proposition}

Analogously to the way the fundamental group is defined, we define the \textbf{uniform fundamental group} $\bm{\check{\pi}_1(X,x_0)}$ to be the group of generalized loops at $x_0$. It is equivalent to the inverse limit of the fundamental groups $\pi_1(R(X,E))$ and, for metric continua, is equivalent to the first shape group \cite[Corollary 6.5]{BDLM}.

\section{Group actions and covering maps}

In this section we address the issue of when a left group action of $G$ on a
uniform space $X$ induces a uniform covering map or a generalized uniform
covering map. As usual, we restrict our attention to faithful group actions
(that means $g\cdot x=x$ for all $x\in X$ implies $g=1\in G$) -
the reason is one can replace $G$ by $G/H$, where $H$ is the stabilizer
of $X$. We do not assume that the action is by uniform equivalences as it is assumed in~\cite{Jam, BP3}.

Recall the action of $G$ on $X$ is \textbf{neutral} \cite[Definition 6.2 on p.90]{Jam}
if for each entourage $E$ of $X$ there is an entourage $F$ of $X$
such that $(x, h\cdot y)\in F$ implies existence of $g\in G$
so that $(g\cdot x,y)\in E$.

\begin{proposition}\label{NeutralActionsAndCLP}
Let $G$ be a group acting on a uniform space $X$.
The projection $p \colon X\to X/G$ has the chain lifting property if and only if
the action is neutral.
\end{proposition}
\proof Suppose the action is neutral and $E$, $F$ are entourages of $X$
such that $(x, h\cdot y)\in F$ implies existence of $g\in G$
so that $(g\cdot x,y)\in E$.
Assume $y_1,\ldots, y_k$
is a $p(F)$-chain in $X/G$. To show it lifts to a $E$-chain in $X$
it suffices to consider the case $k=2$. Suppose $p(x_1)=y_1$.
Choose $(x,y)\in F$ so that $p(x)=y_1$ and $p(y)=y_2$.
There is $g\in G$ satisfying $x=g\cdot x_1$.
Now $(x_1,g^{-1}\cdot y)\in E$ and for $x_2=g^{-1}\cdot y$ one has $p(x_2)=y_2$.
\par Now suppose any $p(F)$-chain in $X/G$ lifts to an $E$-chain.
Given $(x, h\cdot y)\in F$ lift the $p(F)$-chain $\{p(h\cdot y),p(x)\}$
to an $E$-chain $\{y,z\}$. Notice there is $g\in G$ with $z=g\cdot x$.
\endproof

\begin{remark}
Since a map that has chain lifting generates a uniform structure \cite[Proposition 2.8]{BDLM}, if a group $G$ acts neutrally on a uniform space $X$ we give $X/G$ the uniform structure generated by the projection $p \colon X\to X/G$.
\end{remark}

\begin{corollary}\label{UCMAndActions}
The projection $p \colon X\to X/G$ is a uniform covering map if and only if
the action is neutral and there is an entourage $E_0$ of $X$
such that $(x,g\cdot x)\in E_0$ implies $g\cdot x=x$.
\end{corollary}

\begin{definition}\label{PropDiscAction}
Suppose a group $G$ acts on a uniform space $X$.
The action is called {\bf uniformly properly discontinuous} if there exists an entourage $E$ of $X$ such that for all $x\in X$ the inclusion $(x,g\cdot x)\in E$ implies $g=1\in G$.
\end{definition}

Our definition of uniformly properly discontinuous action is weaker than the one used by James~\cite{Jam} since we do not assume that the group acts by uniform equivalences. 

\begin{proposition}\label{properly discontinuous implies UCLP}
Let $G$ be a group acting on a uniform space $X$.
If the action is uniformly properly discontinuous, then the projection $p \colon X\to X/G$ has the uniqueness of chain lifts property.
\end{proposition}
\proof
Consider an entourage $E$ of $X$ such that for all $x\in X$ the inclusion $(x,g\cdot x)\in E$ implies $g=1\in G$. We show that the entourage $E$ is transverse to $p$. Indeed, $p(x)=p(y)$ implies there is $g\in G$ such that $y=g(x)$. Then $(x,y)=(x,g\cdot x)\in E$ implies $g=1\in G$ and $y=x$.
\endproof

\begin{corollary}\label{Neutral+PCIsUCM}
Suppose a group $G$ acts on a uniform space $X$. If the action is neutral and properly discontinuous then the projection $p:X\to X/G$ is a uniform covering map.
\end{corollary}

For the reverse implication we need some continuity of the action. Following Berestovskii and Plaut \cite{BP3}, given an action of $G$ on $X$
and given an entourage $F$ of $X$, we define the set $S_F=\{h\in G \mid (x_h,h\cdot x_h)\in F
\mathrm{\ for \ some\ }x_h\in X\}$.
Let $G_F$ be the subgroup of $G$ generated by $S_F$.

\begin{definition}
Suppose a group $G$ acts on a uniform space $X$. The action is \textbf{small scale uniformly continuous} if for each entourage $E$ of $X$ there is an entourage $F$ of $X$ such that for each $g\in G_F$, $g^{-1}(E)$ is an entourage of $X$.
\end{definition}

Two extreme cases of small scale uniformly continuous actions are:
\begin{itemize}
\item Actions by uniform equivalences --- in this case each map $g$ is uniformly continuous;
\item Properly discontinuous actions (defined below) --- in this case $G_F=1$ for some $F$.
\end{itemize}

\begin{proposition}\label{UCLP implies properly discontinuous}
Let $G$ be a group acting on a chain connected uniform space $X$.
Suppose the action is faithful and small scale uniformly continuous.
If the projection $p \colon X\to X/G$ has the uniqueness of chain lifts property, then
the action is uniformly properly discontinuous.
\end{proposition}
\proof
There is an entourage $E$ transverse to $p$, i.e., $(x,g\cdot x)\in E$ implies $x=g\cdot x$. Since the action is faithful, it is enough to show that if $x=g\cdot x$ for some $x\in X$ and $g\in G$, then $y=g\cdot y$ for any $y\in X$. Assume there is $y\in X$ with $y\ne g\cdot y$. Since $x=g\cdot x$, we have $g\in G_F$ for any entourage $F$. By small scale uniform continuity of the action, there is an entourage $F$ such that $g^{-1}(F)$ is an entourage of $X$ and $F^2\subset E$.
Consider a $F\cap g^{-1}(F)$-chain $\alpha=\{x_0,\dots,x_n\}$ from $x=x_0$ to $y=x_n$. Then $g(\alpha)$ is an $F$-chain. Choose the smallest $i$ satisfying $x_i\ne g(x_i)$. Notice $(x_i,g(x_i))\in F^2$ as $x_{i-1}=g(x_{i-1})$.
Hence $x_i=g(x_i)$, a contradiction.
\endproof

\begin{corollary}
Suppose $G$ is a group acting faithfully on a chain connected uniform space $X$. The action is small scale uniformly continuous and the projection $p:X\to X/G$ is a uniform covering map if and only if the action is neutral and properly discontinuous.
\end{corollary}

Recall $G$ acts on $X$ {\bf equicontinuously} \cite[Definition 6.1 on p.89]{Jam} if for every entourage $E$ of $X$ there is an entourage $F$ of $X$
satisfying $(g\cdot x,g\cdot y)\in E$ for all $g\in G$ and $(x,y)\in F$.
Equivalently, $X$ has a base of entourages $E$ that are $G$-invariant
(that means $(x,y)\in E$ implies $(g\cdot x,g\cdot y)\in E$ for all $g\in G$).
Indeed, $F^\prime=\bigcup\limits_{g\in G} g\cdot F$ is $G$-invariant
and $F^\prime\subset E$ if $(g\cdot x,g\cdot y)\in E$ for all $g\in G$ and $(x,y)\in F$. In \cite{BP3} actions where $X$ has a base of $G$-invariant entourages are called equi-uniform. The definition of discrete actions in \cite{BP3} can be restated as equicontinuous actions that are uniformly properly discontinuous.

\begin{corollary}\label{DiscreteActionsInduceCovers}
A discrete action of $G$ on a uniform space $X$ induces
a uniform covering map $X\to X/G$.
\end{corollary}
\proof
Equi-continuous actions are neutral. Use \ref{Neutral+PCIsUCM}.
\endproof

Now we direct our attention to actions that induce generalized uniform covers. We consider a property that is related to approximate uniqueness of chain lifts as uniformly properly discontinuity is related to uniqueness of chain lifts.

\begin{definition}
Suppose a group $G$ acts on a uniform space $X$. The action has \textbf{small scale bounded orbits} if for each entourage $E$ of $X$ there is an entourage $F$ of $X$ so that the orbits of the induced action of $G_F$ on $X$ are $E$-bounded.
\end{definition}

This property gives a stronger version of small scale continuity.

\begin{definition}
Suppose a group $G$ acts on a uniform space $X$. The action is \textbf{small scale uniformly equicontinuous} if for each entourage $E$ of $X$ there is an entourage $F$ of $X$ such that for each $g\in G_F$, $F\subset g^{-1}(E)$.
\end{definition}

\begin{lemma}\label{SmallOrbitsImpliesEquicontinuity}
Suppose $G$ acts on a uniform space $X$.
If the action has small scale bounded orbits then the action is small scale uniformly equicontinuous.
\end{lemma}
\proof Choose an entourage $E$ such that $E^3\subset D$
and pick $H\subset E$ so that orbits of $G_H$-action are $E$-bounded.
If $(x,y)\in H$ and $g\in G_H$, then $(x,g\cdot x)\in E$
and $(y,g\cdot y)\in E$ therefore $(g\cdot x,g\cdot y)\in E\circ H\circ E\subset D$.
\endproof

Notice that an action is small scale uniformly equicontinuous if and only if for each entourage $E$ of $X$ there is are entourages $F$ and $H$ of $X$ such that for each $g\in G_F$, $H\subset g^{-1}(E)$. Also an action is small scale uniformly continuous if it is small scale uniformly equicontinuous. The following example shows that small scale uniform equicontinuity is indeed weaker than uniform equicontinuity.  It also shows that small scale uniform continuity is weaker than the action being by uniform equivalences.

\begin{example}\label{SmallScaleEquiActionNotEquiUniform}
Let $G$ be the free group generated by $\{x_n\}_{n\ge 1}$
and let $G_k$ be its subgroup generated by $\{x_n\}_{n\ge k}$.
Define $E_k$ as all $(x,y)\in G\times G$ so that $x\cdot y^{-1}\in G_k$
and let $X$ be $G$ with the uniform structure generated by $\{E_k\}_{k\ge 1}$.
\begin{itemize}
\item[a.] The action of $G$ on $X$ defined by left multiplication is small scale uniformly equicontinuous but not equicontinuous.
\item[b.] The projection $X\to X/G$ is a generalized uniform covering map.
\end{itemize}
\end{example}
\proof a. Notice $G_{E_k}=G_k$ and $(x,y)\in E_k$, $g\in G_k$,
imply $(g\cdot x,g\cdot y)\in E_k$ as $(g\cdot x)\cdot (g\cdot y)^{-1}g\cdot (x\cdot y^{-1})\cdot g^{-1}\in G_k$.
\par b. Obviously, $p\colon X\to X/G$ has a chain lifting property.
Suppose $y_1,\ldots ,y_k$ is an $E_n$-chain in $X$.
Thus $g_m=y_m\cdot y_{m+1}^{-1}\in G_n$ for all $m< k$
and $y_1\cdot y_k^{-1}=g_1\cdot\ldots\cdot g_{m-1}\in G_n$.
Hence $(y_1,y_k)\in G_n$ and any two $E_n$-lifts of the constant chain in $X/G$
are $E_n$-close. That means $p$ has approximate uniqueness of chain lifts property.
If $y_1$ and the end-point $y_k$ are fixed then existence of an $E_m$-chain
joining them for any $m$ means $y_1=y_k$. Thus all generalized uniform paths of
$X$ are constant and $p\colon X\to X/G$ is a generalized uniform covering map.
\endproof

\begin{proposition}\label{SubgroupDealImpliesApproximateChainLifting}
Suppose $G$ acts on a uniform space $X$. If the action has small scale bounded orbits then the projection $p\colon X\to X/G$ has approximate uniqueness of
chain lifts.
\end{proposition}
\proof
Given an entourage $E$ of $X$ pick an entourage $F_1$ of $X$
with the property that $(x,g\cdot x)\in E$ for all $x\in X$ and all
elements $g$ of the group $G_{F_1}$.
Also, by~\ref{SmallOrbitsImpliesEquicontinuity} we can assume $(x,y)\in F_1$ implies $(g\cdot x,g\cdot y)\in E$ for all $g\in G_{F_1}$.
Let $D$ be an entourage such that $D^2\subset F_1$.
Choose an entourage $F\subset D$ of $X$ with the property
that $(x,y)\in F$ implies $(g\cdot x,g\cdot y)\in D$ for all $g\in G_D$.

Suppose there are $F$-chains $\alpha=\{x_0,\ldots, x_k\}$ and $\gamma=\{z_0,\ldots,z_k\}$ in $X$ that originate at the same point $x_0=z_0$ and have $x_i=g_i\cdot z_i$ for some $g_i\in G$, $i\leq k$. It suffices to show that $g_k\in G_{F_1}$ since then their endpoints $z_k$ and $x_k$ are $E$-close. Let us use induction on $k$.

Suppose $k=1$. Now $(x_0,x_1)\in F\subset D$ and $(x_0,g_1\cdot x_1)\in F\subset D$ so $(x_1,g_1\cdot x_1)\in D^2\subset F_1$ and $g_1\in G_{F_1}$.

Now suppose $g_{k-1}\in G_{F_1}$. Since $(x_{k-1},x_k)=(g_{k-1}\cdot z_{k-1},x_k)\in F$, $(z_{k-1},g_{k-1}^{-1}\cdot x_k)\in D$. We also have $(z_{k-1},z_k)\in F\subset D$ so $(g_{k-1}^{-1}\cdot x_k,z_k)=(g_{k-1}^{-1}g_k\cdot z_k,z_k)\in D^2\subset F_1$. Therefore $g_{k-1}^{-1}g_k\in G_{F_1}$ so $g_k\in G_{F_1}$.
\endproof

Again, for the reverse implication we need continuity.

\begin{proposition}\label{ApproximateChainLiftingImpliesSubgroupDeal}
Suppose $G$ acts small scale uniformly equicontinuously on a chain connected uniform space $X$.
If $p\colon X\to X/G$ has approximate uniqueness of chain lifts property,
then the action has small scale bounded orbits.
\end{proposition}
\proof
Given an entourage $E$ of $X$ choose an entourage $F_1$ of $X$
such that any two $F_1$-chains $\alpha$ and $\beta$ with a common origin
must be $E$-close if $p(\alpha)=p(\beta)$. Choose an entourage $F$ of $X$ so that $F\subset g^{-1}(E)$ for all $g\in G_F$. 

Suppose $g\in G_F$. We wish to show that $(x,g\cdot x)\in E$ for every $x\in X$. Now $g=h_1\cdots h_n$ where $h_i\in S_F$ for $i\leq n$, say $(y_i,h_i\cdot y_i)\in F$. For every $i\leq n$ we consider an $F$-chain $\alpha_i$ from $\left(\prod_{j=1}^{i-1} h_j\right)\cdot x$ to $y_i$.
Then $(h_i\cdot\alpha_i)^{-1}$ is an $F_1$-chain from $h_i\cdot y_i$ to $\left(\prod_{j=1}^{i} h_j\right)\cdot x$.
Thus, the chain
$$ \alpha=\alpha_1*(h_1\cdot \alpha_1)^{-1}*\dots*\alpha_n*(h_n\cdot\alpha_n)^{-1} $$
is an $F_1$-chain from $x$ to $\left(\prod_{j=1}^{n} h_j\right)\cdot x=g\cdot x$.

Also, for every $i\in\{1,\dots,n\}$ consider the $F_1$-chain $\beta_i=\left(\prod_{j=1}^{i-1} h_j\right)^{-1}\cdot \alpha_i$ starting at $x$.
Notice that $p(\beta_i)=p(\alpha_i)$ and $p(\beta_i^{-1})=p((h_i\cdot \alpha_i)^{-1})$.
Thus, the chain
$$ \beta=\beta_1*(\beta_1)^{-1}*\dots*\beta_n*(\beta_n)^{-1} $$
is an $F$-chain from $x$ to $x$ and $p(\beta)=p(\alpha)$.
Therefore the endpoints of the chains $\alpha$ and $\beta$ are $E$-close.
\endproof

\begin{remark}
We can replace the hypothesis of small scale uniform equicontinuity in \ref{ApproximateChainLiftingImpliesSubgroupDeal} with the action being by uniform equivalences since we only need equicontinuity for the finite set $\left\{ \left( \prod_{j=1}^{i-1} h_j \right)^{\epsilon_i}:i\leq n \mathrm{ , } \epsilon_i=\pm 1 \right\}$. Therefore if $G$ acts on a chain connected uniform space $X$ by uniform equivalences and $p:X\to X/G$ has approximate uniqueness of chain lifts, then the action is small scale uniformly equicontinuous.
\end{remark}

Since the definition of {\bf pro-discrete} actions in \cite{BP3} can be restated as uniformly equicontinuous actions that have small scale bounded orbits, one gets the following:

\begin{corollary}\label{ProDiscreteActionsInduceGenCovers}
If $G$ acts pro-discretely on a chain-connected space $X$, then the induced
map $p\colon X\to X/G$ from $X$ to the orbit space $X/G$ is a generalized uniform covering map if its fibers are complete.
\end{corollary}
\proof By~\ref{SubgroupDealImpliesApproximateChainLifting}
 $p$ has approximate uniqueness of chain lifts property
 and~\ref{GenUnifCovForCompleteFibers} says $p$ is a generalized uniform
 covering map.
\endproof

Recall that an action of a group $G$ on a set $X$ is called {\bf free} if for all $x\in X$ the equality $x=g\cdot x$ implies $g=1\in G$.

\begin{corollary}\label{ActionsInducingGUCMAreFree}
Suppose $G$ acts faithfully on a uniform Hausdorff space $X$.
If the action has small scale bounded orbits,
then the action is free.
\end{corollary}
\proof Suppose $g\cdot x_0=x_0$ for some $g\in G$.
Hence $g\in G_F$ for all entourages $F$, then $(g\cdot x,x)\in E$
for all entourages $E$ and all $x\in X$. Since $X$ is Hausdorff, $g\cdot x=x$
for all $x\in X$ and $g=1$ as the action is faithful.
\endproof

\section{Uniform structures induced by group actions}

There are two well-known uniform
structures on the set $X^G$ of functions from $G$ to $X$:

\begin{enumerate}
\item {\bf Uniform convergence structure} whose base
consists of entourages $E^\ast_X=\{(u,v) | (u(g),v(g))\in E \text{ for all } g\in G\}$,
where $E$ is any entourage of $X$.
\item {\bf Pointwise uniform convergence structure} whose base
consists of entourages $E^\ast_S=\{(u,v) | (u(g),v(g))\in E \text{ for all } g\in S\}$,
where $E$ is any entourage of $X$ and $S$ is any finite subset of $G$.
\end{enumerate}

However, if $G$ acts on a uniform space $X$, one can introduce a third structure.

\begin{definition}\label{SmallScaleUC}
{\bf Small scale uniform convergence structure} on $X^G$
has a base
consisting of entourages $E^\ast=\{(u,v) | (u(g),v(g))\in E \text{ for all } g\in G_E\}$,
where $E$ is any entourage of $X$.
\end{definition}

\begin{proposition}\label{EquiContinuityAndDifferentStructures}
Suppose a group $G$ acts on a uniform space $X$
and $\phi\colon X\to X^G$ is defined by $\phi(x)(g)=g\cdot x$ for $x\in X$ and $g\in G$.
\begin{itemize}
\item[a.] The action is through uniform equivalences if and only if $\phi$ is uniformly continuous
when $X^G$ is equipped with the pointwise convergence structure.
\item[b.] The action is equicontinuous if and only if $\phi$ is uniformly continuous
when $X^G$ is equipped with the uniform convergence structure.
\item[c.] The action is small scale equicontinuous if and only if $\phi$ is uniformly continuous
when $X^G$ is equipped with the small scale uniform convergence structure.
\end{itemize}
\end{proposition}
\proof a. Recall (see \cite{Jam}) $G$ acts through uniform equivalences
if, given $g\in G$, the function $x\to g\cdot x$ from $X$ to $X$ is uniformly
continuous.
\par
Since the projection $X^G\to X$ on the $g$-th coordinate is uniformly
continuous, $\phi$ being uniformly continuous implies that, given $g\in G$,
the map $x\to g\cdot x$ is uniformly continuous.

If the action is through uniform equivalences and $S\subset G$ is finite,
then for any entourage $E$ of $X$ and any $g\in S$ we can find an entourage
$F_g$ of $X$ such that $(x,y)\in F_g$ implies $(g\cdot x,g\cdot y)\in E$.
Pick $F$ such that $F\subset F_g$ for all $g\in S$ and notice
$(x,y)\in F$ implies $(\phi(x),\phi(y))\in E^\ast$, i.e. $\phi$ is uniformly
continuous.
\par b. If the action is equicontinuous and $E$ is an entourage of $X$,
we pick an entourage $F$ with the property $(x,y)\in F$
implies $(g\cdot x,g\cdot y)\in E$ for all $g\in G$. Observe $(x,y)\in F$
implies $(\phi(x),\phi(y))\in E^\ast_X$.
\par Conversely, if $(x,y)\in F$
implies $(\phi(x),\phi(y))\in E^\ast_X$, then $(x,y)\in F$
implies $(g\cdot x,g\cdot y)\in E$ for all $g\in G$ and uniform continuity of $\phi$
implies equicontinuity of the action.
\par c.  If the action is small scale equicontinuous and $E$ is an entourage of $X$,
we pick an entourage $F$ with the property $(x,y)\in F$
implies $(g\cdot x,g\cdot y)\in E$ for all $g\in G_F$. Observe $(x,y)\in F$
implies $(\phi(x),\phi(y))\in E^\ast$.
\par Conversely, if $(x,y)\in F$
implies $(\phi(x),\phi(y))\in E^\ast$, then $(x,y)\in F$
implies $(g\cdot x,g\cdot y)\in E$ for all $g\in G_E$ and uniform continuity of $\phi$
implies small scale equicontinuity of the action.
\endproof

If a group $G$ acts on a uniform space $X$, we consider two ways of creating uniform structure on $G$ (for other uniform structures on $G$ see~\cite[Chapter IV]{Nakano}).

(1) Let $E$ be an entourage of $X$.
Define a symmetric subset $E^*$ of $G\times G$ as follows: $(g,h)\in E^*$ if $(g\cdot x,h\cdot x)\in E$ for all $x\in X$.
Let $\{E^*\}$ be a base of the {\bf uniform convergence structure} on $G$.

(2)
To define a base of {\bf small scale uniform convergence structure} on $G$ we
consider the family $\{\bar E\}$, where $E$ is an entourage of $X$
and $(g,h)\in \bar E$ if $g\cdot h^{-1}\in G_E$.

\begin{corollary}\label{GUCInTermsOfUniformConvergence}
Suppose $G$ acts neutrally and freely on a uniform space $X$
and let $p\colon X\to X/G$ be the induced
map  from $X$ to the orbit space $X/G$.
The following conditions are equivalent:
\begin{itemize}
\item[a.] $p$ is a uniform covering map.
\item[b.] The uniform convergence structure on $G$ is discrete.
\end{itemize}
\end{corollary}
\proof Since the projection $p$ has the chain lifting property
by  \ref{NeutralActionsAndCLP} it suffices to observe
that the uniform convergence structure on $G$ is discrete
if and only if there is an entourage $E$ transverse to $p$.
\endproof

\begin{corollary}\label{GUCMInTermsOfTwoStructures}
Suppose that a group $G$ acts small scale equicontinuously and neutrally on a chain-connected space $X$. Let $p\colon X\to X/G$ be the projection.
Consider the following conditions:
\begin{itemize}
\item[a.] $p$ is a generalized uniform covering map.
\item[b.] The uniform convergence and the small scale uniform convergence structures on $G$ are the same.
\end{itemize}Condition $a)$ implies Condition $b)$.
If the fibers of $p$ are complete, then $b)$ implies $a)$.
\end{corollary}

\proof
a)$\implies$b). By \ref{ApproximateChainLiftingImpliesSubgroupDeal}
given an entourage $E$ of $X$ there is an entourage $F$
such that the orbits of induced action of $G_F$ on $X$
are $E$-bounded. That implies $\bar F\subset E^\ast$.
Indeed, if $g\cdot h^{-1}\in G_F$ and $y=h\cdot x$,
$(g\cdot x,h\cdot x)=((g\cdot h^{-1})\cdot y,y)\in E$.
Since $E^\ast\subset \bar E$ for all $E$, the two structures are the same.

b)$\implies$a). If $\bar F\subset E^\ast$, then the orbits of induced action of $G_F$ on $X$
are $E$-bounded. By \ref{SubgroupDealImpliesApproximateChainLifting}
$p$ has the approximate uniqueness of chain lifts property
and by \ref{GenUnifCovForCompleteFibers} it is a generalized uniform
covering map.
\endproof

\section{Applications to continua theory}
In this section we present a short exposition of a construction
of a homogeneous curve $P$ of J.Prajs \cite{Pra} that is path-connected but not locally connected.
That construction is one of the most interesting results in continua theory in the last decade, so it is worth pursuing a little bit different point of view on it.

The example of $P$ in \cite{Pra} was constructed as the inverse limit of
finite regular covering maps of Menger sponge $M$ over itself. By \cite{Lab1} it leads to a pro-discrete action on $P$
and as such $P$ is theoretically realizable
in the form of the orbit space $GP(M,m_0)/K$ for some subgroup $K\subset
\check\pi_1(M,m_0)$ (see \cite{Lab2}).
The purpose of this section is to construct such group $K$
and show the desired properties of $GP(M,m_0)/K$ using the theory developed
in this paper.
Focusing attention on the group $K$ allows to isolate its features that are responsible
for certain properties of $P$. For example, our proof of path-connectedness of $P$
is a real simplification in comparison to \cite{Pra}.
\par Another simplification is the use of Bestvina's theory \cite{Bes}
of Menger manifolds, most notably the result that $UV^0$-maps
between $\mu^1$-manifolds are near-homeomorphisms (consult \cite{Anc}
and \cite{Bro} for basic results on inverse sequences and near-homeomorphisms).
That particular theorem serves very well as the lighthouse for what we are doing
here regarding the Menger sponge $\mu^1=M$.

First, it makes sense to derive general properties of $GP(M,m_0)/K$
depending on $K$.

\begin{theorem}\label{MainHomogeneousThm}
Suppose $M$ is a locally connected continuum,
$\mathcal{H}$ is a transitive group of homeomorphisms
 of $M$, and $K$ is a closed subgroup of
 $\check\pi_1(M,m_0)$.
\begin{itemize}
\item[a.] If $\check\pi_1(M,m_0)/K$ is compact, then $M_K=GP(M,m_0)/K$ is compact.
\item[b.] If  the set of
free generalized loops generated by $K$ is $\mathcal{H}$-invariant, then $M_K$ is homogeneous.
\item[c.] If $\check\pi_1(M,m_0)/K$ is Abelian, then $M_K$ is path-connected.
\item[d.] If there is a closed neighborhood $N$ of $m_0$
such that the image of the composition $\check\pi_1(N,m_0)\to \check\pi_1(M,m_0)\to
\check\pi_1(M,m_0)/K$
has empty interior, then $M_K$ is not locally path-connected.
\end{itemize}
\end{theorem}
\proof
Let $q\colon GP(M,m_0)\to M_K$ and $p\colon M_K\to M$ be the quotient maps.
\par {\bf a)}.
Since $GP(M,m_0)$ is metrizable, the space $M_K$ has a countable base of entourages and we might as well
use metric reasoning in that case.
It suffices to show $p$ is a closed map (in the topological category).
If $p$ is not closed, there is a sequence $x_n\in M_K$
forming a closed subset of $M_K$
such that $p(x_n)$ converges to $p(x_0)$ but $p(x_n)\ne p(x_0)$
for all $n > 0$. Express $x_n$ as $q(\alpha_n)$
and find paths $\beta_n$ from $p(x_n)$ to $p(x_0)$ whose diameter
tends to $0$ as $n\to\infty$.
Put $y_n=q(\alpha_n\ast\beta_n)$ for $n > 0$ and notice $p(y_n)=p(x_0)$.
Without loss of generality we may assume
$y_n$ converges to $y_0$ as $p^{-1}(p(x_0))$ is compact.
Since diameters of $\beta_n$ approach $0$,
$y_0$ must be the limit of $x_n$ as well. Hence $y_0=x_n$ for some $n > 0$
and $p(x_0)=p(x_n)$, a contradiction.
\par
{\bf b)}. By a free generalized loop generated by $K$ we mean any loop
of the form $\alpha^{-1}\ast \beta\ast \alpha$, where $\beta\in K$
and $\alpha$ is a generalized path joining $m_0$ to some $x$.
The set $T$ of free generalized loops being $\mathcal{H}$-invariant means
$\tilde h(T)\subset T$ for any $h\in \mathcal{H}$.
\par
Since $\check\pi_1(M,m_0)/K$ serves as the group of deck transformations
of $p$, all we need to show is that for any generalized path $\alpha$ from $m_0$
to some $x\ne m_0$ there is a homeomorphism $H$ of $M_K$
so that $H(q(m_0))=q(\alpha)$.
Choose a homeomorphism $h$ of $M$ sending $m_0$ to $x$.
Define $H$ by $H(q(\beta))=q(\alpha\ast \tilde h(\beta))$.
$H$ is well-defined exactly when $T$ is $\mathcal{H}$-invariant.
Indeed, if $\gamma\in K$,
then $\alpha\ast \tilde h(\gamma\ast \beta)\alpha\ast \tilde h(\gamma)\ast \tilde h(\beta)\alpha\ast \tilde h(\gamma)\ast \alpha^{-1}\ast \alpha\ast \tilde h(\beta)$
and $q(\alpha\ast \tilde h(\gamma\ast \beta))=q(\alpha\ast \tilde h(\beta))$ as
$\alpha\ast \tilde h(\gamma)\ast \alpha^{-1}\in K$.
\par {\bf c)}.
Let us show first that the elements of the fiber $p^{-1}(m_0)$
can be joined by a path in $M_K$ to $q(m_0)$.
Such element corresponds to an element $\gamma\in \check\pi_1(M,m_0)/K$.
Since that group is Abelian,
the homomorphism $\check\pi_1(M,m_0)\to \check\pi_1(M,m_0)/K$
factorizes through the abelianization of $\check\pi_1(M,m_0)$,
the first \v Cech homology group $\check H_1(M)$ of $M$
(see \cite{Dyd3}). As $\check H_1(M)$ is the image of the singular
first homology group $H_1(M)$ of $M$ (see \cite{EdaKaw})
and $H_1(M)$ is the abelianization of $\pi_1(M,m_0)$,
we may realize $\gamma$ as a real loop in $M$.
That loop lifts to a path in $M_K$ starting at $q(m_0)$
an ending in the desired element of $p^{-1}(m_0)$.
\par If $x\in M\setminus\{m_0\}$, we choose a path $\lambda$ from $x$ to $m_0$
in $M$. Given an element $y$ of $p^{-1}(x)$ we can lift $\lambda$
to a path in $M_K$ starting at $y$ an ending in $p^{-1}(m_0)$.
That path can be continued to $q(m_0)$.
That completes the proof of $M_K$ being path-connected.

\par {\bf d)}.
If $M_K$ is locally path connected, there is a path connected neighborhood
$U$ of $q(m_0)$ in $p^{-1}(N)=\check\pi_1(M,m_0)/K$.
All the paths from $q(m_0)$ to elements of $p^{-1}(N)\cap U$
generate loops in $(N,m_0)$ whose images via
the composition $\check\pi_1(N,m_0)\to \check\pi_1(M,m_0)\to
\check\pi_1(M,m_0)/K$ form exactly $(\check\pi_1(M,m_0)/K)\cap U$,
a contradiction.
\endproof

\subsection{Z-sets and $\frac{1}{2}$Z-circles}
From now on let $M$ be the Menger sponge.
Recall a {\bf Z-set} in $M$ (see \cite{Bes}) is a closed subset of $M$
such that the identity map $id_M$ belongs to the closure
of the set of maps $f\colon M\to M$ missing $A$.
Notice any near-homeomorphism $h\colon M\to M$
sends a Z-set of the form $h^{-1}(A)$ to a Z-set.

A subset $C$ of $M$ homeomorphic to the circle
 is called a {\bf $\frac{1}{2}$Z-circle} if it has a {\bf Z-arc-resolution}:
 there is a copy $M_1$ of the Menger sponge and an arc $A_1$ that is a
 Z-set in $M_1$ such that the pair $(M,C)$ is obtained from $(M_1,A_1)$
 by gluing two disjoint non-degenerate intervals $B_1$ and $C_1$ of $A_1$.
 If $q_1\colon M_1\to M$ is the quotient map corresponding to gluing,
 notice every sub-arc of $C$ missing $q_1(B_1)$
is a Z-set in $M$ and none of sub-arcs of $q_1(B_1)$ is a Z-set in $M$.
Thus every $\frac{1}{2}$Z-circle $C$ has a specific structure
in the form of the arc $q_1(B_1)$ and the closure of its complement.
Any homeomorphism between two $\frac{1}{2}$Z-circles
preserving that structure will be called a {\bf $\frac{1}{2}$Z-homeomorphism}.

\par The general strategy from now on is to follow the well-known scheme of things
for Z-sets (homeomorphism extension theorems) and apply
it to $\frac{1}{2}$Z-circles via Z-arc-resolutions.

 For a $\frac{1}{2}$Z-circle define $\kappa_C\in H^1(M;\Z/2)$ as follows:
 Pick a retraction $r_1$ of $M_1$ onto $A_1$ and notice it induces
a retraction $s_1\colon M\to C$.
The composition $g\circ s_1\colon M\to RP^\infty$,
where $g\colon C\to RP^\infty$ generates $H^1(C;\Z/2)$,
is the desired $\kappa_C$.

\begin{lemma}\label{CharClassIsIndependent}
$\kappa_C$ does not depend neither on the choice of gluing nor retraction $r_1$.
\end{lemma}
\proof
If $r^{\prime}_1\colon M_1\to A_1$ is another retraction, it is homotopic
to $r_1$ rel. $A_1$, so the induced retraction $s^{\prime}_1\colon M\to C$
is homotopic to $s_1$.
\par Suppose another set of objects $M_2$, $A_2$, $B_2$, $C_2$, and $r_2$ is given
with the corresponding gluing denoted by $q_2\colon M_2\to M$.
Notice $q_2(B_2)=q_1(B_1)$ as every sub-arc of $C$ missing $q_2(B_2)$
is a Z-set in $M$ and none of sub-arcs of $q_2(B_2)$ is a Z-set in $M$.
Pick a homeomorphism $h\colon B_1\to B_2$ so that $h\circ q_2=q_1$,
extend it over $C_1$ to $h\colon C_1\cup B_1\to C_2\cup B_2$,
and then to $h\colon A_1\to A_2$. Finally, extend it to $h\colon M_1\to M_2$.
Given a retraction $r_1\colon M_1\to A_1$ we can produce
the corresponding $r_2\colon M_2\to A_2$ as $h\circ r_1\circ h^{-1}$
and notice $r_1$ and $r_2$ induce the same retraction $M\to C$.
 \endproof

\begin{corollary}\label{HomeoAndCharClass}
Suppose $C$ is a $\frac{1}{2}$Z-circle in $M$.
If $h\colon M\to M$ is a homeomorphism, then
$h(C)$ is a $\frac{1}{2}$Z-circle in $M$ and
$\kappa_{h(C)}\circ h=\kappa_C$.
\end{corollary}

The proof of \ref{CharClassIsIndependent} yields another useful fact
in conjunction with the standard extension of homeomorphism between
Z-sets in $M$ (see \cite{Bes}).

\begin{corollary}\label{HalfZHomeoExtension}
A homeomorphism between two $\frac{1}{2}$Z-circles in $M$
extends to a homeomorphism of $M$ onto $M$ if and only if
it is a $\frac{1}{2}$Z-homeomorphism.
\end{corollary}

With a bit more work one gets the following:
\begin{corollary}\label{ZandHalfZHomeoExtensionPlusDenseSets}
Suppose $g\colon C_1\to C_2$ is a $\frac{1}{2}$Z-homeomorphism between two $\frac{1}{2}$Z-circles in $M$ and $h\colon A_1\to A_2$ is a homeomorphism
of Z-sets in $M$, $A_1\subset M\setminus C_1$ and $A_2\subset M\setminus C_2$.
Given two countable dense sets $D_1$ and $D_2$ in $M$,
$D_1\subset M\setminus (C_1\cup A_1)$ and $D_2\subset M\setminus (C_2\cup A_2)$,
there is a homeomorphism $H\colon M\to M$
extending both $g$ and $h$ such that $H(D_1)=D_2$.
\end{corollary}
\proof By passing to Z-arc-resolutions of both $C_1$ and $C_2$ one reduces
\ref{ZandHalfZHomeoExtensionPlusDenseSets} to the case
of Z-sets only in which case one uses local pushes (see the proofs later on)
utilizing the fact $M$ has a base of open sets $U$ whose boundary
$\partial(U)$ is a Z-set in $U\cup\partial(U)$.
\endproof

\subsection{The construction of the group $K$}

What we need is a large set of linearly independent elements $\kappa_C$ of $H^1(M;\Z/2)$.
Notice that $\{\kappa_{C(i)}\}_{i=1}^n$ are linearly independent if curves $C(i)$ are mutually disjoint.
Indeed, $\alpha=\sum\limits_{i=1}^n \kappa_{C(i)}$ cannot be $0$ as $\alpha | C(i)$ is
the generator of $H^1(C(i);\Z/2)$ for each $i$.

As in \cite{Pra} we choose a countable family of Z-arcs $\{A_i\}_{i=1}^\infty$ in $M$
that forms a null-set (that means only finitely many of them have diameter bigger than a given $\epsilon$ for any $\epsilon > 0$), such that $\bigcup\limits_{i=1}^\infty A_i$
is dense in $M$ and we use those arcs to form a quotient map
$q\colon M\to M_1$ by gluing two sub-arcs of each $A_i$ in order to get
a $\frac{1}{2}$Z-circle $C_i$. In comparison to \cite{Pra} the technique
of near-homeomorphisms is very handy here: using Bestvina's characterization
\cite{Bes} of $M$ as a locally connected continuum that has disjoint arc property,
one can easily see that performing finally many gluing operations yields a Menger sponge again.
For infinite set of gluing operations one forms $M_n$ as the result of
contracting arcs $A_i$, $i > n$, to points and then gluing parts of arcs $A_i$, $i\leq n$.
One has a natural projection $M_{n+1}\to M_n$ that is $UV^0$, hence a near-homeomorphism.
Consequently, the inverse limit of those projections (which is $M_1)$
is homeomorphic to $M$.

\begin{definition}
Given a family $\mathcal{F}=\{C(s)\}_{s\in S}$ of
mutually disjoint $\frac{1}{2}$Z-circles of $M$ define the following
objects:
\begin{enumerate}
\item $L(\mathcal{F})$ is the set of all loops $f\colon S^1\to M$ such that $\kappa_{C(s)}\circ f$
is null-homotopic for all $s\in S$.
\item $\mathcal{H}(\mathcal{F})$ is the group of all homeomorphisms $h$ so that $h(L(\mathcal{F}))\subset L(\mathcal{F})$
(by $h(L(\mathcal{F}))$ we mean loops of the form $h\circ f$, $f\in L(\mathcal{F}))$.
\item The group $K(\mathcal{F})$ consists of all generalized loops in $\check\pi_1(M,m_0)$
whose image in $\check H_1(M)$
is represented by a loop $f\in L(\mathcal{F})$.
\item $T(\mathcal{F})$ is the set of free generalized loops of $K(\mathcal{F})$.
\end{enumerate}
\end{definition}

Notice that elements of $T(\mathcal{F})$, after being sent to $\check H_1(M)$,
are represented by exactly the same loops $f\in L(\mathcal{F})$, so $T(\mathcal{F})$ being $\mathcal{H(\mathcal{F})}$-invariant
is easy.

\begin{lemma}
If $\mathcal{F}$ forms a null-set,
then $K(\mathcal{F})$ is closed and $\check\pi_1(M,m_0)/K(\mathcal{F})$ is compact
Abelian. Moreover, if $\bigcup\limits_{s\in S}C(s)$ is dense in $M$,
then for any proper closed subset $N$ of $M$ containing $m_0$
the image of the composition $\check\pi_1(N,m_0)\to \check\pi_1(M,m_0)\to
\check\pi_1(M,m_0)/K(\mathcal{F})$
has empty interior.
\end{lemma}
\proof
Let $\alpha\colon M\to RP_S=\prod\limits_{s\in S} RP^2$
be the diagonal map induced by $\kappa_{C(s)}$, $s\in S$.
Notice $K(\mathcal{F})$ is exactly the kernel of $\check\pi_1(\alpha)\colon \check\pi_1(M,m_0)\to \check\pi_1(RP_S,c_0)=\check H_1(RP_S)$, so it suffices to show $\check H_1(\alpha)\colon \check H_1(M)\to \check H_1(RP_S)$ being epimorphic. That follows from continuity of $\check H_1(\alpha)$
and the fact finite unions of images of loops $C(s)$ form a dense subset of $\check H_1(RP_S)$.
\par
If $N$ is proper, then infinitely many sets $C(t)$ are contained in $M\setminus U$,
$U$ a neighborhood of $N$. For any such $t\in S$ the image
of $\check\pi_1(N,m_0)\to \check\pi_1(M,m_0)\to
\check\pi_1(M,m_0)/K(\mathcal{F})=\prod\limits_{s\in S}(\Z/2)_s$
has $t$-coordinate equal $0$, hence cannot have non-empty interior
in $\prod\limits_{s\in S}(\Z/2)_s$.
\endproof

\begin{lemma}\label{TransitivityLemma}
Suppose $\mathcal{F}=\{C(s)\}_{s\in S}$ forms a null-set
of
mutually disjoint $\frac{1}{2}$Z-circles
and $\bigcup\limits_{s\in S}C(s)$ is dense in $M$.
$\mathcal{H}(\mathcal{F})$ acts on $M$ transitively if
any decomposition of $M$ into some members of  $\mathcal{F}$
and singletons yields $M$ up to homeomorphism.
 \end{lemma}
\proof
First observe any homeomorphism $h$ that leaves $\bigcup\limits_{s\in S}C(s)$
invariant belongs to $\mathcal{H}(\mathcal{F})$
(use \ref{HomeoAndCharClass}), so \ref{TransitivityLemma} follows from
the following two facts:
\par {\bf Fact 1}: Any two points $x$ and $y$ not belonging
to $\bigcup\limits_{s\in S}C(s)$ can be connected by such homeomorphism
that is an identity on a given Z-set $A$ of $M$ contained in $M\setminus \bigcup\limits_{s\in S}C(s)\setminus \{x,y\}$
(see an outline later on).
\par {\bf Fact 2}: If $x\in C(s_0)$, then one can find a homeomorphism $h$
arbitrarily close to $id_M$ such that $h(\bigcup\limits_{s\ne s_0}C(s))= \bigcup\limits_{s\ne s_0}C(s)$
and $h(x)\notin \bigcup\limits_{s\in S}C(s)$ (see an outline later on).
\par Indeed, by Fact 2 one can make $h$ so close to $id_M$ that $\kappa_{C(s_0)}$ and $\kappa_{C(s_0)}\circ h$
are homotopic resulting in $h(L(\mathcal{F}))\subset L(\mathcal{F})$.
\par Fact 2 is actually a consequence of Fact 1 as follows.
Let $D$ be the decomposition of $M$ into $C(s)$, $s\ne s_0$, and singletons
with the resulting quotient map $q\colon M\to M/D$.
Since $M/D$ is homeomorphic to $M$, one can find a closed
neighborhood $B\approx M$ in $M/D$ of $q(x)$ whose boundary $\partial B$
is a Z-set in $B$ and is contained in $q(M\setminus \bigcup\limits_{s\ne s_0}C(s)$.
Now $A=q^{-1}(\partial B)$ is a $Z$-set in $q^{-1}(B)\approx M$, so using Fact 1
we can find a homeomorphism $h$ of $M$ fixing $M\setminus int(q^{-1}(B))$ permuting
those $C(s)$ that are contained in $q^{-1}(B)$. That homeomorphism can be made as close
as possible to $id_M$ by adjusting the size of $B$.
\par Fact 1 can be reduced to:
\par {\bf Fact 3}: Suppose $A$ and $B$ are $\frac{1}{2}$Z-circles in $M$
 and $h_1\colon A\to B$ is a $\frac{1}{2}$Z-homeomorphisms.
 Let $M_1=M/A$ and $M_2=M/B$.
Given any homeomorphism $g\colon M_1\to M_2$
sending $A/A$ to $B/B$ and given any neighborhood $U$ of $B/B$
there is a homeomorphism
$h\colon M\to M$ coinciding with $g$ outside of $g^{-1}(U)$
so that $h | A=h_1$.
Moreover, if $D_A$ is a countable dense set in $M\setminus A$
and $D_B$ is a countable dense set in $M\setminus B$ and $g(D_A)=D_B$,
we can accomplish $h(D_A)=D_B$.

\par Let us show {\bf Fact 3 $\implies$ Fact 1}.
Let $M_1$ be the decomposition of $M$ into $C(s)$, $s\in S$, and singletons.
Let $D_1$ be the image of $\bigcup\limits_{s\in S} C_s$ in $M_1$ via the quotient map
$q_1\colon M\to M_1$. $D_1$ is a countable set, so there is
a homeomorphism $h_1\colon M_1\to N_1=M_1$ satisfying $h_1(q_1(x))=q_1(y)$,
$h_1 | q_1(A)=id$,
and $h_1(D_1)=D_1$.
Assume $S$ is the set of natural numbers and let $M_n$ be the quotient
of $M$ via the decomposition of $M$ into $C(s)$, $s\ge n$, and singletons.
Let $q_n\colon M\to M_n$ be the quotient map and let $D_n=q_n(\bigcup\limits_{s\ge n}C(s))$.
Let $q^{n+1}_n\colon M_{n+1}\to M_n$ be the natural quotient map.
We want to create homeomorphisms $h_n\colon M_n\to N_n$ for $n\ge 2$,
where $N_n$ is the decomposition of $M$ into $C(t_1),\ldots,C(t_n)$
and singletons. We have analogous quotient maps $p_n\colon M\to N_n$
and $p^{n+1}_n\colon N_{n+1}\to N_n$, and analogous countable dense sets $E_n$
in $N_n$, so we require
 $h_n(D_n)=E_n$ and $p^n_{n-1}\circ h_n$ is $\epsilon_n$-close
to $h_{n-1}\circ q^n_{n-1}$, where $\{\epsilon_n\}_{n\ge 1}$ is
going to be constructed to ensure convergence of $h_n$ to a homeomorphism
$h\colon M\to M$ (see \cite{Anc} or \cite{Bro}).
Moreover, we want $h_n(q_n(\bigcup\limits_{s\in S}C(s)))= p_n(\bigcup\limits_{s\in S}C(s))$
and $h(q_n(x))=q_n(y)$.
\par $h_n$ is obtained by picking a small (of diameter less than $\frac{\epsilon_{n-1}}{2}$) neighborhood $B\approx M$ of the point
$h_{n-1}(q_{n-1}(C(n)))$ so that the boundary $\partial B$ of $B$ is a Z-set in $B$
and does not intersect $p_{n-1}(\bigcup\limits_{s\in S}C(s))$.
Let $t_n$ be the element of $S$ satisfying $h_{n-1}(q_{n-1}(C(n)))=p_{n-1}(C(t_n))$.
Put $B^{\prime}=h_{n-1}^{-1}(B)$.
On the set $(q^n_{n-1})^{-1}(B^{\prime})$ we construct a homeomorphism
sending

\begin{enumerate}
\item the boundary of $(q^n_{n-1})^{-1}(B^{\prime})$ to the boundary
of $(p^n_{n-1})^{-1}(B)$ as determined by $h_{n-1}$,

\item $q_n(C(n))$ onto $p_n(C(t_n)$,

\item $D_n\cap (q^n_{n-1})^{-1}(B^{\prime})$ onto $E_n\cap (p^n_{n-1})^{-1}(B)$,

\end{enumerate}

and then we patch it with the lift of $h_{n-1}$.

\endproof

\end{document}